\documentstyle[12pt,amssymb]{amsart}

\numberwithin{equation}{section}

\newtheorem{lem}{Lemma}[section]

\newtheorem{thm}[lem]{Theorem}

\begin{document}

\title[Self-dual Koornwinder-Macdonald polynomials]{Self-dual
Koornwinder-Macdonald polynomials}
\author[J. F. van Diejen]{J. F. van Diejen}
\address{$\!$Department of Mathematical Sciences, University of Tokyo,
Hongo~7-3-1, Bunkyo-ku, Tokyo~113, Japan}

\subjclass{Primary 33D45, 05E05, 05E35; Secondary 05A19}
\keywords{basic orthogonal polynomials in several variables,
duality properties, recurrence relations, normalization constants,
generalized Selberg-type integrals}

\thanks{Work supported by the Japan Society for the Promotion of Science
(JSPS).}

\date{(July 1995)}

\maketitle

\begin{abstract}
We prove certain duality properties and present recurrence relations for a
four-parameter family of self-dual Koornwinder-Macdonald polynomials. The
recurrence relations are used to verify Macdonald's normalization conjectures
for these polynomials.
\end{abstract}

\section{Introduction}\label{sec1}
In a to date unpublished but well-known manuscript, Macdonald introduced
certain families of multivariable orthogonal polynomials associated with
(admissible pairs of integral) root systems and conjectured the values of the
normalization constants turning these polynomials into an orthonormal system
\cite{mac:orthogonal}.
Recently, Cherednik succeeded in verifying Macdonald's normalization
conjectures in the case of reduced root systems
(and admissible pairs of the form $(R,R^\vee)$) using a technique involving
so-called shift operators \cite{che:double}.
Previously, this same technique had enabled Opdam to prove the normalization
conjectures for a degenerate case ($q\rightarrow 1$) of the Macdonald
polynomials known as the Heckman-Opdam-Jacobi polynomials
\cite{opd:some,hec:elementary}.

Meanwhile, a generalization of Macdonald's construction for the nonreduced root
system $BC_n$---resulting in a multivariable version of the famous Askey-Wilson
polynomials \cite{ask-wil:some}---was presented by Koornwinder
\cite{koo:askey}.
It turns out that {\em all}
Macdonald polynomials associated with classical (i.e.,
non-exceptional) root systems may be seen as special cases of these
multivariable Askey-Wilson polynomials \cite[Sec. 5]{die:commuting}
(type $A$ by picking the highest-degree homogeneous parts of the polynomials
and types $B$, $C$, $D$, and $BC$, by specialization of the parameters).

In the present paper, we will prove certain duality properties and recurrence
relations for (a four-parameter subfamily of) the Koornwinder-Macdonald
multivariable Askey-Wilson polynomials, which enable one to verify the
corresponding Macdonald conjectures for the (ortho)normalization constants
also in this (more general) situation.
Our approach does not involve shift operators but rather exploits the fact that
the polynomials are joint eigenfunctions of a family of commuting difference
operators that was introduced by the author in Ref.~\cite{die:commuting} (see
also Ref.~\cite{die:diagonalization}). By duality, these difference operators
give rise to a system of recurrence relations from which, in turn, the
normalization constants follow.

The same method employed here was used already several years ago by Koornwinder
when verifying similar duality properties and normalization constants for the
Macdonald polynomials related to the root system $A_{n}$
\cite{koo:self,mac:symmetric}.
(In this special case, though, the validity of the normalization conjectures
had also been checked by Macdonald himself.) The $A_n$-type Macdonald
polynomials constitute a multivariable generalization of the $q$-ultraspherical
polynomials \cite{ask-wil:some} (to which they reduce for $n=1$). The present
paper may thus be regarded as an extension of Koornwinder's methods in
Ref.~\cite{koo:self} (see also Ref.~\cite[Ch. 6]{mac:symmetric}) to the
multivariable Askey-Wilson level, or, if one prefers, as an extension
from type $A$ root systems to type $BC$ root systems.

\section{Koornwinder-Macdonald polynomials}\label{sec2}
The Koornwinder-Macdonald multivariable Askey-Wilson polynomials are
characterized by a weight function of the form
\begin{equation}\label{weight}
\Delta (x)\; =\;
\prod\begin{Sb} 1\leq j<j^\prime \leq n \\
                \varepsilon ,\varepsilon^\prime =\pm 1 \end{Sb}
d_v (\varepsilon\, x_j+\varepsilon^\prime\, x_{j^\prime})\;
\prod\begin{Sb} 1\leq j\leq n \\ \varepsilon =\pm 1 \end{Sb}
d_w (\varepsilon\, x_j),
\end{equation}
where
\begin{gather*}
d_v(z)\; =\;
\frac{(e^{-\alpha z}; q)_\infty}
     {( q^g e^{-\alpha z}; q)_\infty},\;\;\;\;\;\;\;\;\;\;\;\;\;\;
     q=e^{-\alpha\beta},   \\[1ex]
   d_w(z) =
\frac{(e^{-\alpha z},\,
      -e^{-\alpha z},\,
       q^{1/2} e^{-\alpha z},\,
      -q^{1/2} e^{-\alpha z};\,
                                    q)_\infty    }
     {(q^{g_0} e^{-\alpha z},\,
      -q^{g_1} e^{-\alpha z},\,
       q^{(g_2+1/2)}e^{-\alpha z},\,
      -q^{(g3+1/2)} e^{-\alpha z};\,
                                    q)_\infty    }  ,
\end{gather*}
and the $q$-shifted factorials are defined, as usual, by
$(a;q)_\infty=\prod_{l=0}^\infty (1-a q^l)$ and
$(a_1,\ldots ,a_k;q)_\infty = (a_1;q)_\infty \cdots (a_k;q)_\infty$.
To ensure the convergence of the infinite products contained in $\Delta $
\eqref{weight} it will be assumed that $\alpha ,\beta > 0$ (so $1<q<1$); in
addition, we will also assume $g,g_r \geq 0$, $r=0,1,2,3$.

Let $\{ m_\lambda (x)\}_{\lambda \in \Lambda}$ denote the basis consisting of
even and permutation symmetric exponential
monomials (or monomial symmetric functions)
\begin{equation}
m_\lambda (x) \; =\sum_{\lambda^\prime \in W \lambda}
e^{\alpha \sum_{j=1}^n \lambda_j^\prime x_j},
\;\;\;\;\;\; \lambda \in \Lambda = \{ \lambda \in {\Bbb Z}^n \; |\; \lambda_1
\geq \lambda_2 \geq \cdots \geq \lambda_n \geq 0\; \} , \label{monomial}
\end{equation}
with $W$ being the group generated by permutations and sign flips
of $x_j$, $j=1,\ldots ,n$ ($W\cong S_n \ltimes ({\Bbb Z}_2)^n$).
The monomial basis can be partially ordered by defining
for $\lambda ,\lambda^\prime \in \Lambda$ \eqref{monomial}
\begin{equation}\label{ord}
\lambda^\prime \leq \lambda \;\;\;\;\;\;\;\; \text{iff} \;\;\;\;\;\;\;\;
\sum_{1\leq j\leq m} \lambda_j^\prime \leq \sum_{1\leq j\leq m} \lambda_j
\;\;\;\;\; \text{for} \;\;\;\;\; m=1,\ldots ,n
\end{equation}
(and $\lambda^\prime < \lambda$ iff $\lambda^\prime \leq \lambda$ and
$\lambda^\prime \neq \lambda$).
The Koornwinder-Macdonald polynomials $p_\lambda(x)$,
$\lambda\in\Lambda$ can now be introduced as the (unique) trigonometric
polynomials satisfying
\begin{itemize}
\item[i.] $\displaystyle p_\lambda (x) = m_\lambda (x) +
\sum_{\lambda^\prime \in \Lambda, \lambda^\prime < \lambda }\;
 c_{\lambda ,\lambda^\prime }\: m_{\lambda^\prime}(x)$,
\ \ \ \ \ $\displaystyle c_{\lambda ,\lambda^\prime}\in {\Bbb C}$;
\item[ii.]
$\displaystyle  \langle p_\lambda  ,  m_{\lambda^\prime}\rangle_\Delta  =0$
\ \ if \ \  $\displaystyle \lambda^\prime < \lambda$,
\end{itemize}
where $\langle\cdot ,\cdot\rangle_\Delta$ denotes the inner product
determined by
\begin{equation}\label{ip}
\langle m_\lambda ,m_{\lambda^\prime} \rangle_\Delta =
\left( \frac{\alpha}{2\pi} \right)^n
\int_{-\pi/\alpha}^{\pi/\alpha}\!\!\!\!\!\!\cdots
               \int_{-\pi/\alpha}^{\pi/\alpha}
m_\lambda (ix)\, \overline{m_{\lambda^\prime}(ix)}\,  \Delta (ix)\,
dx_1\cdots dx_n
\end{equation}
(and extended by bilinearity).
In other words, the polynomial $p_\lambda (ix)$ consists of the monomial
$m_\lambda (ix)$ minus its orthogonal projection in
$L^2( ]-\frac{\pi}{\alpha}, \frac{\pi}{\alpha} ]^n, \Delta (ix)dx_1\cdots
dx_n)$ onto the finite-dimensional subspace
span$\{ m_{\lambda^\prime}(ix)\}_{\lambda^\prime\in \Lambda,\,
\lambda^\prime < \lambda}$.

It is of course possible to extend this construction of multivariable
polynomials determined by Conditions i. and ii. to a more general class of
weight functions than the one considered here. In general, however, the
resulting polynomials will not be orthogonal (except for $n=1$) because the
ordering in Eq.~\eqref{ord} is not a total ordering (unless $n=1$). (A priori
the construction only guarantees
that $p_\lambda (x)$ and $p_{\lambda^\prime}(x)$ be orthogonal if $\lambda$ and
$\lambda^\prime$ are comparable with respect to the ordering in
Eq.~\eqref{ord}.)
Still, it turns out \cite{koo:askey} that for the weight function
$\Delta $ \eqref{weight} the
corresponding polynomials indeed do constitute an orthogonal system for
arbitrary $n$:

{\em Orthogonality}
\begin{equation}\label{ort}
\langle p_\lambda ,p_{\lambda^\prime} \rangle_\Delta  =0\;\;\;\;\;\;\;\;\;\;
\text{if}\;\;\;\;\;\; \lambda \neq \lambda^\prime .
\end{equation}
This feature should be looked upon as a very restrictive property of the weight
function $\Delta $ \eqref{weight}.

{\em Remark:} The relation between our parameters and the parameters employed
by Koornwinder reads (cf.~\cite[Eqs. (5.1), (5.2)]{koo:askey})
\begin{equation}\label{parel}
t=q^g,\;\;\;\;\; a=q^{g_0},\;\;\;\;\; b=-q^{g_1},\;\;\;\;\;
c=q^{(g_2+1/2)},\;\;\;\;\; d=-q^{(g_3+1/2)} .
\end{equation}
Furthermore, Koornwinder fixes the period the trigonometric functions
to be $2\pi (i)$, i.e., he puts $\alpha =1$.

\section{Difference equations}\label{sec3}
Another special property of the polynomials associated with $\Delta $
\eqref{weight} is that they satisfy a second order difference equation
\cite{koo:askey}. (This difference equation is in fact instrumental in the
orthogonality proof.) It can be written as
\begin{equation}\label{diffeq1}
\begin{split}
\begin{align}
\sum_{1\leq j\leq n} \Bigl(
w(x_j)\prod_{k\neq j} v(x_j +x_k)\, v(x_j -x_k) \:
                    &[p_\lambda (x+\beta e_j) - p_\lambda (x)]+ \\ \nonumber
w(-x_j)\prod_{k\neq j} v(-x_j +x_k)\, v(-x_j -x_k) \:
                    &[p_\lambda (x-\beta e_j) - p_\lambda (x)] \Bigr)
\end{align} \\
 = E(\rho +\lambda  )\: p_\lambda (x) ,
\end{split}
\end{equation}
with
\begin{gather*}
v(z)= \frac{\text{sh} \frac{\alpha}{2} (\beta g + z)}
                 {\text{sh} (\frac{\alpha}{2} z)},    \\
w(z)= \frac{\text{sh} \frac{\alpha}{2} (\beta g_0 + z)}
           {\text{sh} (\frac{\alpha}{2} z)}
       \frac{\text{ch} \frac{\alpha}{2} (\beta g_1 + z)}
            {\text{ch} (\frac{\alpha}{2} z)}
      \frac{\text{sh} \frac{\alpha}{2} (\beta g_2 +\frac{\beta}{2} + z)}
           {\text{sh} \frac{\alpha}{2}(\frac{\beta}{2} + z)}
      \frac{\text{ch} \frac{\alpha}{2} (\beta g_3 +\frac{\beta}{2} + z)}
           {\text{ch} \frac{\alpha}{2}(\frac{\beta}{2} + z)}, \\
 \intertext{and}
 E(y)=2\sum_{1\leq j\leq n} \left(
\text{ch} (\alpha\beta y_j)-\text{ch} (\alpha\beta \rho_j) \right) ,\\
\rho = \sum_{1\leq j\leq n} \rho_j\, e_j,
\;\;\;\;\;\;\;\;\;\; \rho_j = (n-j)\, g + (g_0+g_1+g_2+g_3)/2.
\end{gather*}
(The vector $e_j$ denotes the $j$-th unit element of standard basis in
${\Bbb R}^n$.) For $n=1$ this difference equation reduces to the well-known
difference equation for the Askey-Wilson polynomials \cite{ask-wil:some}.

In Ref.~\cite{die:commuting} it was shown that for arbitrary
number of variables $n$ the above
difference equation can be extended to a system of $n$ independent difference
equations of order $2r$, $r=1,\ldots ,n$, respectively. This system is
explicitly given by:

{\em Difference equations}
\begin{equation}\label{diffeqr}
\sum\begin{Sb} J\subset \{ 1,\ldots ,n\} ,\, 0\leq |J|\leq r \\
               \varepsilon_j=\pm 1,\; j\in J      \end{Sb}\!\!\!\!
U_{J^c,\, r-|J|}(x)\,  V_{\varepsilon J,\, J^c}(x)\,
p_\lambda (x+\beta e_{\varepsilon J})= E_r (\rho +\lambda  )\: p_\lambda (x),
\;\;\;\; r=1,\ldots , n,
\end{equation}
with
\begin{gather*}
\begin{align*}
V_{\varepsilon J,\, K}(x) =
\prod_{j\in J} w(\varepsilon_jx_j) &
\prod\begin{Sb} j,j^\prime \in J \\ j<j^\prime \end{Sb}
v(\varepsilon_jx_j+\varepsilon_{j^\prime}x_{j^\prime})
 v(\varepsilon_jx_j+\varepsilon_{j^\prime}x_{j^\prime}+\beta )\label{VJ}\\
 \times & \prod\begin{Sb} j\in J \\ k\in K \end{Sb} v(\varepsilon_j x_j+x_k)
v(\varepsilon_j x_j -x_k), \nonumber
\end{align*} \\
\begin{align*}
U_{K,p}(x)= (-1)^p \!\! \sum\begin{Sb} L\subset K,\, |L|=p \\
                                  \varepsilon_l =\pm 1,\; l\in L \end{Sb}\;
\prod_{l\in L} w(\varepsilon_l x_l)\,
& \prod\begin{Sb} l,l^\prime \in L \\ l<l^\prime \end{Sb}
v(\varepsilon_lx_l+\varepsilon_{l^\prime}x_{l^\prime})
v(-\varepsilon_lx_l-\varepsilon_{l^\prime}x_{l^\prime}-\beta ) \label{UKm}\\
\times &
\prod\begin{Sb} l\in L \\ k\in K\setminus L \end{Sb} v(\varepsilon_l x_l+x_k)
v(\varepsilon_l x_l -x_k) ,\nonumber
\end{align*} \\
\intertext{and}
 E_r\: (y) =
2^r\!\!\!\sum\begin{Sb} J\subset \{ 1,\ldots ,n\} \\ 0\leq |J|\leq r\end{Sb}
\!\! (-1)^{r-|J|} \Bigl(  \prod_{j\in J} \text{ch}(\alpha\beta y_j)
\!\!\!\! \sum_{r\leq l_1\leq\cdots\leq l_{r-|J|}\leq n}\!\!\!\!
\text{ch}(\alpha\beta \rho_{l_1})\cdots \text{ch}(\alpha\beta
\rho_{l_{r-|J|}}) \Bigr) .
\end{gather*}
In the above formulas $|J|$ represents the number of elements of
$J\subset \{ 1,\ldots ,n \}$ and
\begin{equation}
e_{\varepsilon J}=\sum_{j\in J} \varepsilon_j e_j\;\;\;\;\;\;\;\;\;\;\;\;
\;\;\; (\varepsilon_j =\pm 1).
\end{equation}
Furthermore, we used the conventions that empty products are equal to one,
$U_{K,\, p}=1$ if $p=0$, and the second sum in $E_r(y)$ is equal to
one if $|J|=r$.
For $r=1$, the difference equation in Eq.~\eqref{diffeqr} reduces to
that of Eq.~\eqref{diffeq1}.

{\em Remarks: i.} Equation~\eqref{diffeqr} may be interpreted as a system of
eigenvalue equations
\begin{equation}
(D_r\, p_\lambda )(x) = E_r (\rho +\lambda )\:  p_\lambda (x) ,
\end{equation}
for a family of $n$ independent commuting difference operators of the form
\begin{equation}\label{ados}
D_r(x)=\sum\begin{Sb} J\subset \{ 1,\ldots ,n\} ,\, 0\leq |J|\leq r \\
               \varepsilon_j=\pm 1,\; j\in J      \end{Sb}\!\!\!\!
U_{J^c,\, r-|J|}(x)\,  V_{\varepsilon J,\, J^c}(x)\,
T_{\varepsilon J,\, \beta},
\;\;\;\; r=1,\ldots , n,
\end{equation}
with $T_{\varepsilon J,\, \beta}=\prod_{j\in J} T_{\varepsilon_j j,\beta}$
and
\begin{equation*}
(T_{\pm j,\beta} f)(x_1,\ldots ,x_n)=
f(x_1,\ldots ,x_{j-1},x_j\pm\beta ,x_{j+1},\ldots ,x_n).
\end{equation*}
The commutativity of $D_1(x),\ldots ,D_n(x)$ follows \cite{die:commuting}
from the fact that
the difference operators are simultaneously diagonalized by the basis
$\{ p_\lambda (x) \}_{\lambda \in \Lambda}$.
Two other properties of the operators that follow immediately from their
diagonalization by the Koornwinder-Macdonald polynomials are the triangularity
with respect to the partially ordered monomial basis
$\{ m_\lambda (x)\}_{\lambda \in \Lambda}$, and the symmetry with respect to
the
inner product $\langle \cdot ,\cdot \rangle_\Delta$ \eqref{ip} (the eigenvalues
$E_r (\rho +\lambda)$ are real).
More precisely, one has
({\em Triangularity})
\begin{equation}\label{tria}
(D_r\, m_\lambda ) (x) = E_r(\rho +\lambda )\; m_\lambda (x) +
\sum_{\lambda^\prime \in \Lambda , \lambda^\prime < \lambda }
[D_r]_{\lambda ,\lambda^\prime } m_{\lambda^\prime}(x),\;\;\;\;
[D_r]_{\lambda ,\lambda^\prime } \in {\Bbb C},
\end{equation}
and
({\em Symmetry})
\begin{equation} \label{sym}
 \langle D_r\, m_\lambda , m_{\lambda^\prime} \rangle_\Delta =
\langle m_\lambda , D_r\, m_{\lambda^\prime}\rangle_\Delta.
\end{equation}

{\em ii.}  For $r=1$ the operator $D_r(x)$ \eqref{ados} reduces to
(cf. the l.h.s. of Eq.~\eqref{diffeq1})
\begin{equation*}
\begin{align}\label{dop1}
D(x)= \sum_{1\leq j\leq n} \Bigl(
w(x_j)\prod_{k\neq j} v(x_j +x_k)\, v(x_j -x_k) \:
                    &(T_{j,\beta} -1)+ \\ \nonumber
w(-x_j)\prod_{k\neq j} v(-x_j +x_k)\, v(-x_j -x_k) \:
                    &(T_{-j,\beta}-1) \Bigr)  .
\end{align}
\end{equation*}
With the aid of this operator the following useful representation for the
Koornwinder-Macdonald polynomials can be given (cf. Ref.~\cite{mac:orthogonal})
\begin{equation}\label{urep}
p_\lambda (x) = \left(
\prod_{\lambda^\prime \in \Lambda , \lambda^\prime < \lambda }
\frac{ D(x)-E(\rho +\lambda^\prime )}{E(\rho +\lambda )-
E(\rho +\lambda^\prime )}
\right) m_\lambda (x) .
\end{equation}
(Using the Triangularity \eqref{tria} and Symmetry \eqref{sym} of
$D(x)$~\eqref{dop1}, it is not difficult to verify that the r.h.s. of Eq.
\eqref{urep} satisfies the Conditions i. and ii. in Section~\ref{sec2}.
Notice also that the denominators in the r.h.s. of Eq.~\eqref{urep} do not
vanish
because $E(\rho +\lambda^\prime) < E(\rho +\lambda)$ if
$\lambda^\prime < \lambda $ (cf. Ref.~\cite[Lemma 5.1]{die:commuting}).)
It follows from this representation for $p_\lambda (x)$ that the
Koornwinder-Macdonald polynomials are rational in
$\exp (\alpha\beta g)$ and $\exp (\alpha\beta g_r)$, $r=0,1,2,3$.
Hence, they may be extended uniquely to nonnegative
(or even complex) values of the parameters $g,g_r$.
In view of the analytic dependence on the parameters it is clear that the
resulting polynomials then satisfy the Difference equations \eqref{diffeqr}
for all the values of $g, g_r$, $r=0,1,2,3$.

\section{Duality and recurrence relations}\label{sec4}
{\em Note:} In this section, we will drop the condition that
the parameters $g, g_r$, $r=0,1,2 ,3$ be nonnegative (cf. Remark {\em ii}. of
the previous section).

In order to describe the duality relations it is convenient to introduce
certain dual polynomials $p^*_\lambda (x)$, $\lambda\in \Lambda$. These dual
polynomials are again Koornwinder-Macdonald polynomials but with a slightly
different parametrization.
Specifically, the parameters of $p^*_\lambda (x)$ are related to those of
$p_\lambda (x)$ by
\begin{equation}\label{repa}
\begin{array}{c}
\alpha^*=\beta ,\\ \beta^* =\alpha ,\\ g^* =g,
\end{array}
\;\;\;\;\;\;\;\;\;\;
\left(
\begin{matrix} g^*_0 \\ g^*_1 \\ g^*_2 \\ g^*_3 \end{matrix}
                                                                   \right)
= \frac{1}{2} \left(
\begin{array}{rrrr}
1  &  1  &  1  &  1  \\
1  &  1  &  -1  &  -1  \\
1  &  -1  &  1  &  -1  \\
1  &  -1  &  -1  &  1
\end{array}             \right)
\left( \begin{matrix} g_0 \\ g_1 \\ g_2 \\ g_3 \end{matrix} \right) .
\end{equation}
Notice that the reparametrization in Eq.~\eqref{repa} is involutive, i.e.,
$p^{**}_\lambda (x) =p_\lambda (x)$.

Furthermore, instead of working with monic polynomials we go over to a
different normalization by introducing
\begin{equation}\label{reno}
\tilde{p}_\lambda(x) =\frac{p_\lambda (\beta x)}
                   {p_\lambda (\beta\rho^* )},
\;\;\;\;\;\;\;\;\;\;\;\;\;\;\;
\tilde{p}^*_\lambda(x) =\frac{p^*_\lambda (\alpha x)}
                         {p^*_\lambda (\alpha \rho )}
\end{equation}
($\rho^* = \sum_{j=1}^n \rho_j^* e_j$ with
$\rho_j^*=(n-j)g^* +(g_0^*+\cdots +g_3^*)/2$), where we have also rescaled the
arguments of the trigonometric polynomials such that $\tilde{p}_\lambda(x)$ and
$\tilde{p}^*_\lambda (x)$ have the same (imaginary) period
in $x_j$, $j=1,\ldots ,n$ (viz. $2\pi i/(\alpha \beta )$).
Of course, the renormalization in Eq.~\eqref{reno} only makes sense provided
$p_\lambda (\beta\rho^* )$ and $p^*_\lambda (\alpha \rho )$ do not vanish. This
is guaranteed, at least for generic parameters, by the following lemma.
\begin{lem}\label{nonzero} For generic parameters one has
$$p_\lambda (\beta\rho^* ), p^*_\lambda (\alpha \rho )\neq 0 .$$
\end{lem}
\begin{pf}
For $g,g_0,\ldots ,g_3=0$, the polynomial $p_\lambda (x)$ reduces
to the monomial symmetric function $m_\lambda (x)$ and $\rho^* =0$.
Thus, it is clear that for this special choice of the parameters
$p_\lambda (\beta \rho^*) \neq 0$. But then the same follows for
$g,g_0,\ldots ,g_3$ in an open dense subset of ${\Bbb R}$ (or ${\Bbb C}$)
because
of the analytic dependence of $p_\lambda (x)$ on the parameters
(cf. Remark {\em ii.} in Section~\ref{sec3}).
The analogous statement for $p^*_\lambda(\alpha \rho)$ follows by duality.
\end{pf}
In Section~\ref{sec5} the value of $p_\lambda (\beta\rho^*)$ will be computed
explicitly. We will then see that $p_\lambda (\beta\rho^*)$ is positive for all
nonnegative values of parameters $g, g_r$, $r=0,1,2,3$.

The matrix in Eq.~\eqref{repa} relating $(g^*_0,\ldots ,g^*_3)^t$ and
$(g_0,\ldots ,g_3)^t$ has eigenvalues $+1$ (with multiplicity three) and $-1$
(with multiplicity one). The invariant subspace corresponding to the eigenvalue
$+1$ consists of the hyperplane $g_0-g_1-g_2-g_3=0$.
For parameters in this hyperplane one has $g^*_r=g_r$ ($r=0,1,2 ,3$)
and $\tilde{p}^*_\lambda (x) =\tilde{p}_\lambda (x)$. In other words,
for these parameters the polynomials $\tilde{p}_\lambda (x)$ are self-dual.
In the rest of the paper we will always assume that the {\em Self-duality
condition}
\begin{equation}\label{self-dual}
g_0-g_1-g_2-g_3=0
\end{equation}
is satisfied (unless explicitly stated otherwise).

After these preparations we are now ready to formulate the duality theorem,
which relates the value of $\tilde{p}_\lambda (x)$ in the point $\rho^*+\mu$ to
value of $\tilde{p}^*_\mu (x)$ in the point $\rho +\lambda$
($\lambda  , \mu \in \Lambda$).
\begin{thm}[duality relations]\label{dthm}
Let $\lambda ,\mu \in \Lambda$~\eqref{monomial}. Then the renormalized
Koorn\-winder-Macdonald
polynomials $\tilde{p}_\lambda (x)$ and $\tilde{p}^*_\mu (x)$ \eqref{reno}
satisfy the relation
\begin{equation}\label{dr}
 \tilde{p}_\lambda (\rho^*+\mu) = \tilde{p}^*_\mu (\rho +\lambda ).
\end{equation}
\end{thm}
Theorem~\ref{dthm} was conjectured by Macdonald in Ref.~\cite{mac:some}
(without imposing the Self-duality condition~\eqref{self-dual}).
In the present self-dual set-up it is of course not necessary to distinguish
between the polynomials $\tilde{p}_\lambda (x)$ and the dual polynomials
$\tilde{p}^*_\lambda(x)$ as both polynomials coincide when
Condition~\eqref{self-dual} holds.
However, we have chosen to keep this distinction in our notation
because it is expected that with the present formulation all results
remain valid also when the Self-duality condition~\eqref{self-dual}
is not satisfied (cf. Remark~\ref{sub2} of Section~\ref{sec7}).

Before going to the proof of Theorem~\ref{dthm}, which is relegated to
Section~\ref{sec6}, let us first discuss some important consequences
of these duality relations.
The main point is that the Difference equations \eqref{diffeqr} together with
the Duality relations~\eqref{dr} imply a system of recurrence relations for the
Koornwinder-Macdonald polynomials. To see this one first substitutes $x=\alpha
(\rho +\lambda)$ in the difference equations for the dual polynomial
$p^*_\mu (x)$, $\mu\in\Lambda$.
After dividing by $p_\mu^*(\alpha\rho)$ and invoking of
Definition~\eqref{reno} one arrives at an equation of the form
\begin{equation}\label{dualeqs}
E^*_r  ( \rho^* +\mu )\: \tilde{p}^*_\mu (\rho +\lambda )=\!\!\!\!\!\!\!
\sum\begin{Sb} J\subset \{ 1,\ldots ,n\} ,\, 0\leq|J|\leq r \\
               \varepsilon_j=\pm 1,\; j\in J      \end{Sb}\!\!\!\!
\tilde{U}_{J^c,\, r-|J|}(\rho +\lambda)\,
\tilde{V}_{\varepsilon J,\, J^c}(\rho+\lambda)\,
\tilde{p}^*_\mu (\rho+\lambda+ e_{\varepsilon J}) ,
\end{equation}
where
\begin{gather*}
\begin{align*}
\tilde{V}_{\varepsilon J,\, K}(x) =
\prod_{j\in J} \tilde{w}(\varepsilon_jx_j) &
\prod\begin{Sb} j,j^\prime \in J \\ j<j^\prime \end{Sb}
\tilde{v}(\varepsilon_jx_j+\varepsilon_{j^\prime}x_{j^\prime})
\tilde{v}(\varepsilon_jx_j+\varepsilon_{j^\prime}x_{j^\prime}+1 ) \\
 \times & \prod\begin{Sb} j\in J \\ k\in K \end{Sb}
\tilde{v}(\varepsilon_j x_j+x_k)
\tilde{v}(\varepsilon_j x_j -x_k), \nonumber
\end{align*} \\
\begin{align*}
\tilde{U}_{K,p}(x)= (-1)^p \!\! \sum\begin{Sb} L\subset K,\, |L|=p \\
                                  \varepsilon_l =\pm 1,\; l\in L \end{Sb}\;
\prod_{l\in L} \tilde{w}(\varepsilon_l x_l)\,
& \prod\begin{Sb} l,l^\prime \in L \\ l<l^\prime \end{Sb}
\tilde{v}(\varepsilon_lx_l+\varepsilon_{l^\prime}x_{l^\prime})
\tilde{v}(-\varepsilon_lx_l-\varepsilon_{l^\prime}x_{l^\prime}-1 ) \\
\times &
\prod\begin{Sb} l\in L \\ k\in K\setminus L \end{Sb}
\tilde{v}(\varepsilon_l x_l+x_k)
\tilde{v}(\varepsilon_l x_l -x_k) ,\nonumber
\end{align*}
\end{gather*}
with
\begin{gather*}
\tilde{v}(z)= \frac{\text{sh} \frac{\alpha\beta }{2} (g^* + z)}
                 {\text{sh} (\frac{\alpha\beta }{2} z)},    \\
\tilde{w} (z)= \frac{\text{sh} \frac{\alpha\beta}{2} (g^*_0 + z)}
           {\text{sh} (\frac{\alpha\beta}{2} z)}
       \frac{\text{ch} \frac{\alpha\beta}{2} (g^*_1 + z)}
            {\text{ch} (\frac{\alpha\beta}{2} z)}
      \frac{\text{sh} \frac{\alpha\beta}{2} (g^*_2 +\frac{1}{2} + z)}
           {\text{sh} \frac{\alpha\beta}{2}(\frac{1}{2} + z)}
      \frac{\text{ch} \frac{\alpha\beta}{2} (g^*_3 +\frac{1}{2} + z)}
           {\text{ch} \frac{\alpha\beta}{2}(\frac{1}{2} + z)}
\end{gather*}
(and $E^*_r(y)$ is the dual of $E_r(y)$ in Eq.~\eqref{diffeqr}, i.e., with
$\rho_j$ replaced by $\rho^*_j$).
One may restrict the summation in Eq.~\eqref{dualeqs} to those index sets
$J\subset \{ 1,\ldots ,n\}$ and configurations of signs $\varepsilon_j$,
$j\in J$ for which $\lambda +e_{\varepsilon J}\in
\Lambda$~\eqref{monomial} because of the following lemma.
\begin{lem}\label{res} Let $\lambda \in \Lambda =
\{ \lambda \in {\Bbb Z}^n \; |\; \lambda_1
\geq \lambda_2 \geq \cdots \geq \lambda_n \geq 0\; \}$.
For generic parameters one has
\begin{equation*}
 \tilde{V}_{\varepsilon J,\, J^c}(\rho+\lambda ) = 0
\;\;\;\;\;\; \text{iff}\;\;\;\;\; \lambda +e_{\varepsilon J} \not\in \Lambda .
\end{equation*}
\end{lem}
\begin{pf}
Suppose $\tilde{V}_{\varepsilon J,\, J^c}(\rho+\lambda ) = 0$. Then we are in
one of the following situations.
\begin{enumerate}
\item $\tilde{w}(\varepsilon_j (\rho_j+\lambda_j))=0$ for some $j\in J$.

For generic parameters this only happens when $\lambda_j=0$, $\varepsilon_j=-1$
(and $j=n$). So, $\lambda +e_{\varepsilon J}\not\in \Lambda$ because
$\lambda_j+\varepsilon_j <0$.

\item $\tilde{v}(\varepsilon_j(\rho_j+\lambda_j)+
\varepsilon_{j^\prime}(\rho_{j^\prime}+\lambda_{j^\prime}))=0$
for some $j,j^\prime \in J$ with $j<j^\prime$.

For generic parameters this only happens when
$-\varepsilon_j=\varepsilon_{j^\prime}=1$, $\lambda_j=\lambda_{j^\prime}$
(and $j^\prime =j+1$).
So, $\lambda +e_{\varepsilon J}\not\in \Lambda$ because
$\lambda_j+\varepsilon_j < \lambda_{j^\prime}+\varepsilon_{j^\prime}$
with $j<j^\prime$.

\item $\tilde{v}(\varepsilon_j(\rho_j+\lambda_j)+
\varepsilon_{j^\prime}(\rho_{j^\prime}+\lambda_{j^\prime})+1)=0$
for some $j,j^\prime \in J$ with $j<j^\prime$.

For generic parameters this only happens when
$-\varepsilon_j=\varepsilon_{j^\prime}=1$, $\lambda_j=\lambda_{j^\prime}+1$
(and $j^\prime =j+1$).
So, $\lambda +e_{\varepsilon J}\not\in \Lambda$ because
$\lambda_j+\varepsilon_j < \lambda_{j^\prime}+\varepsilon_{j^\prime}$
with $j<j^\prime$.

\item $\tilde{v}(\varepsilon_j(\rho_j+\lambda_j)+
\varepsilon_{k}(\rho_{k}+\lambda_{k}))=0$
for some $j\in J$ and $k\not\in J$.

For generic parameters this only happens when $\varepsilon_j=-\varepsilon_k=1$,
$\lambda_j=\lambda_k$ and $j=k+1$, or
$\varepsilon_j=-\varepsilon_k=-1$, $\lambda_j=\lambda_k$ and $j=k-1$.
So, $\lambda +e_{\varepsilon J}\not\in \Lambda$ because
$\lambda_j+\varepsilon_j < \lambda_k$ with $j\in J$, $k\not\in J$ and $j<k$, or
$\lambda_j+\varepsilon_j > \lambda_k$ with $j\in J$, $k\not\in J$ and $j>k$.
\end{enumerate}

Conversely, if $\lambda^\prime=\lambda +e_{\varepsilon J}\not\in \Lambda$ then
either
$\lambda_n^\prime < 0$ or there exist a $j\in \{ 1,\ldots ,n-1\}$ such that
$\lambda_j^\prime<\lambda_{j+1}^\prime$. The first case, i.e.
$\lambda_n^\prime< 0$, can occur only if
$n\in J$ and $\lambda_n=0$, $\varepsilon_n =-1$. The vanishing of
$\tilde{V}_{\varepsilon J,\, J^c}(\rho+\lambda )$ then follows from
the vanishing of $\tilde{w}(\varepsilon_n (\rho_n+\lambda_n))$.
In the second case, i.e. $\lambda_j^\prime < \lambda_{j+1}^\prime$, it is
convenient to distinguish the following three situations.
\begin{enumerate}
\item $j,j+1\in J$.

One has $\lambda_j^\prime
=\lambda_j+\varepsilon_j<\lambda_{j+1}+\varepsilon_{j+1}=\lambda_{j+1}^\prime$
only if $0\leq \lambda_j-\lambda_{j+1}\leq 1$ and
$-\varepsilon_j=\varepsilon_{j+1}=1$. The vanishing of $\tilde{V}_{\varepsilon
J,\, J^c}(\rho+\lambda )$ then follows because either
$\tilde{v}(\varepsilon_j(\rho_j+\lambda_j)+
\varepsilon_{j^\prime}(\rho_{j^\prime}+\lambda_{j^\prime}))=0$
or
$\tilde{v}(\varepsilon_j(\rho_j+\lambda_j)+
\varepsilon_{j^\prime}(\rho_{j^\prime}+\lambda_{j^\prime})+1)=0$
(depending on whether $\lambda_j=\lambda_{j+1}$ or
$\lambda_j =\lambda_{j+1}+1$).

\item $j\in J$, $j+1\not\in J$.

One has
$\lambda_j^\prime=\lambda_j+\varepsilon_j<\lambda_{j+1}=\lambda_{j+1}^\prime$
only if
$\lambda_j=\lambda_{j+1}$ and $\varepsilon_j=-1$. The vanishing of
$\tilde{V}_{\varepsilon J,\, J^c}(\rho+\lambda )$ then follows because
$\tilde{v}(\varepsilon_j(\rho_j+\lambda_j)+
(\rho_{j^\prime}+\lambda_{j^\prime}))=0$

\item $j\not\in J$, $j+1\in J$.

One has $\lambda^\prime_j=\lambda_j <
\lambda_{j+1}+\varepsilon_{j+1}=\lambda_{j+1}^\prime$ only if
$\lambda_j=\lambda_{j+1}$ and $\varepsilon_{j+1}=1$. The vanishing of
$\tilde{V}_{\varepsilon J,\, J^c}(\rho+\lambda )$ then follows because
$\tilde{v}(-(\rho_j+\lambda_j)+
\varepsilon_{j^\prime}(\rho_{j^\prime}+\lambda_{j^\prime}))=0$.
\end{enumerate}
(The a priori fourth situation $j,j+1\not\in J$ does not occur
because in that case
$\lambda^\prime_j=\lambda_j \geq \lambda_{j+1}=\lambda^\prime_{j+1}$.)
\end{pf}
In order to restrict the sum in Eq.~\eqref{dualeqs} to the index sets and signs
with $\lambda +e_{\varepsilon J}\in \Lambda$ it is of course sufficient to know
that $\tilde{V}_{\varepsilon J,\, J^c}(\rho+\lambda ) = 0$ if
$\lambda +e_{\varepsilon J}\not\in \Lambda$. It is clear from the proof of
Lemma~\ref{res} that this is actually true for all values of the parameters
(the genericity of the parameters was needed only when proving the converse
statement).

After restricting the sum, we may apply the duality theorem
(Theorem~\ref{dthm}) to obtain the equation
\begin{equation}\label{rec1}
E^*_r (\rho^*+\mu )\: \tilde{p}_{\lambda} (\rho^* +\mu )=\!\!\!\!\!\!\!
\sum\begin{Sb} J\subset \{ 1,\ldots ,n\} ,\, 0\leq|J|\leq r \\
               \varepsilon_j=\pm 1,\; j\in J;\;
                \lambda +e_{\varepsilon J} \in \Lambda  \end{Sb}\!\!\!\!
\tilde{U}_{J^c,\, r-|J|}(\rho +\lambda)\,
\tilde{V}_{\varepsilon J,\, J^c}(\rho +\lambda)\,
\tilde{p}_{\lambda +e_{\varepsilon J}} (\rho^* +\mu ) .
\end{equation}
Equation~\eqref{rec1} describes an equality between trigonometric polynomials
that is satisfied for all $\rho^*+\mu$, $\mu\in \Lambda$.
But then the equality must hold identically for all values of the argument and
we arrive at the following theorem.
\begin{thm}[recurrence relations]\label{recthm}
The renormalized Koornwinder-Macdonald polynomials $\tilde{p}_\lambda (x)$
\eqref{reno}, $\lambda \in \Lambda$~\eqref{monomial},
satisfy a system of recurrence relations of the form
\begin{equation}\label{recr}
E^*_r (x)\: \tilde{p}_{\lambda} (x)=\!\!\!\!\!\!\!
\sum\begin{Sb} J\subset \{ 1,\ldots ,n\} ,\, 0\leq|J|\leq r \\
               \varepsilon_j=\pm 1,\; j\in J;\;
                \lambda +e_{\varepsilon J} \in \Lambda  \end{Sb}\!\!\!\!
\tilde{U}_{J^c,\, r-|J|}(\rho +\lambda)\,
\tilde{V}_{\varepsilon J,\, J^c}(\rho +\lambda)\,
\tilde{p}_{\lambda +e_{\varepsilon J}} (x) ,
\end{equation}
$r=1,\ldots ,n$.
\end{thm}

For $r=1$ the Recurrence relation \eqref{recr} becomes
\begin{equation*}
\begin{align*}
&E^*(x)\: \tilde{p}_\lambda (x) = \\
&\sum_{\begin{Sb} 1\leq j\leq n \\ \lambda +e_j \in \Lambda \end{Sb}}
\!\!\tilde{w}(\rho_j\!+\!\lambda_j)\prod_{k\neq j}
\tilde{v}(\rho_j\!+\!\lambda_j\!+\!\rho_k\!+\!\lambda_k)\,
\tilde{v}(\rho_j\!+\!\lambda_j\! -\!\rho_k\!-\!\lambda_k) \:
                    [\tilde{p}_{\lambda +e_j}(x) - \tilde{p}_\lambda (x)] + \\
\nonumber
&\sum_{\begin{Sb} 1\leq j\leq n \\ \lambda -e_j \in \Lambda \end{Sb}}
\!\! \tilde{w}(-\rho_j\!-\!\lambda_j)\prod_{k\neq j}
\tilde{v}(-\rho_j\!-\!\lambda_j\! +\!\rho_k\!+\!\lambda_k)\,
\tilde{v}(-\rho_j\!-\!\lambda_j\! -\!\rho_k\!-\!\lambda_k) \:
                    [\tilde{p}_{\lambda -e_j} (x) - \tilde{p}_\lambda (x)]
\end{align*}
\end{equation*}
($E^*(y)=E^*_1(y)$).
In this special case the recurrence formula was conjectured by Macdonald
\cite{mac:some}; after specialization to one variable ($n=1$) this formula
reduces to the three-term recurrence relation for the Askey-Wilson polynomials
\cite{ask-wil:some}.

In the one-variable case both duality relations and self-duality
can be easily checked directly through the explicit representation
of $\tilde{p}_l(x)$ in terms of the basic hypergeometric ${}_4\phi_3$-series
(cf. \cite[Eq. (5.8)]{ask-wil:some} and Relation~\eqref{parel}):
\begin{equation}\label{bhrep}
\tilde{p}_l(x) = {}_4\phi_3
\left(
\begin{matrix}
q^{-l},\; q^{2g_0^* +l},\; q^{g_0-x},\; q^{g_0+x} \\ [0.5ex]
-q^{g_0+g_1},\; q^{g_0+g_2+1/2},\; -q^{g_0+g_3+1/2}
\end{matrix} \; ; q,q \right) .
\end{equation}
It is clear from Representation~\eqref{bhrep} that the duality relations
$\tilde{p}_l(g_0+m)=\tilde{p}_m^*(g_0^*+l)$ hold actually
without restriction on the parameters (recall that for $n=1$
one has $\rho =g_0^*$ and $\rho^* =g_0$) and that $\tilde{p}_l(x)$ is self-dual
iff $g_0^*=g_0$, or equivalently, iff $g_0-g_1-g_2-g_3=0$
(notice that $g_0+g_r=g_0^*+g_r^*$ for $r=1,2,3$).

{\em Remark:} Recurrence relations of the type as in Theorem~\ref{recthm} are
sometimes called generalized Pieri formulas (after similar formulas for
the Schur functions) \cite{mac:symmetric}. Similarly, the duality relations of
Theorem~\ref{dthm} are also referred to as symmetry relations.

\section{Normalization}\label{sec5}
In this section we exploit the recurrence relations of Theorem~\ref{recthm}
to evaluate the normalization constants $p_\lambda (\beta \rho^*)$ and
$\langle p_\lambda ,p_\lambda \rangle_\Delta $ (giving rise to polynomials
satisfying the Duality relations \eqref{dr} and orthonormal polynomials,
respectively).
First some notation is needed. Let
\begin{equation}\label{Dp1}
\tilde{\Delta}_+(x) = \prod_{1\leq j<j^\prime \leq n}
\tilde{d}_v(x_j+x_{j^\prime})\, \tilde{d}_v(x_j-x_{j^\prime})
\prod_{1\leq j\leq n} \tilde{d}_w(x_j),
\end{equation}
with
\begin{gather*}
\tilde{d}_v(z)\; =\; q^{-g^* z/2}
\frac{(q^z; q)_\infty}
     {(q^{(g^*+z)}; q)_\infty}\;\;\;\;\;\;\;\;\;\;\;\;\;\;\;\;\;
(q=e^{-\alpha\beta}),\\[1ex]
\tilde{d}_w(z) = q^{-(g^*_0+\cdots+g^*_3)z/2}
\frac {(q^z,\,
      -q^z,\,
       q^{(1/2+z)},\,
      -q^{(1/2+z)};\, q)_\infty    }
     {(q^{(g^*_0+z)},\,
      -q^{(g^*_1+z)},\,
       q^{(g^*_2+1/2+z)},\,
      -q^{(g^*_3+1/2+z)};\, q)_\infty    }  ,
\end{gather*}
and let
\begin{equation}\label{Dp2}
\hat{\Delta}_+(x) = \prod_{1\leq j<j^\prime \leq n}
\hat{d}_v(x_j+x_{j^\prime})\, \hat{d}_v(x_j-x_{j^\prime})
\prod_{1\leq j\leq n} \hat{d}_w(x_j),
\end{equation}
with
\begin{gather*}
\hat{d}_v(z)\; =\; q^{ g^* (z+1)/2}
\frac{(q^{(z+1)};\, q)_\infty}
     {(q^{(-g^*+z+1)};\, q)_\infty},\\[1ex]
%\begin{align*}
\hat{d}_w(z) = q^{(g_0^*+\cdots+g^*_3)(z+1)/2}
\frac {(q^{(z+1)},\,
      -q^{(z+1)},\,
       q^{(3/2+z)},\,
      -q^{(3/2+z)};\, q)_\infty    }
     {(q^{(-g^*_0+z+1)},\,
      -q^{(-g^*_1+z+1)},\,
       q^{(-g^*_2+3/2+z)},\,
      -q^{(-g^*_3+3/2+z)};\, q)_\infty    }  .
%\end{align*}
\end{gather*}
(The function $\hat{\Delta}_+(x)$ \eqref{Dp2} is obtained from
$\tilde{\Delta}_+(x)$ \eqref{Dp1} by substituting $g^*\rightarrow -g^*$,
$g^*_r\rightarrow -g^*_r$ ($r=0,1,2 ,3$), and $z\rightarrow z+1$.)

{}From the relations $\tilde{d}_v(z+1)=\tilde{v}(z)\tilde{d}_v(z)$,
$\tilde{d}_w(z+1)=\tilde{w}(z)\tilde{d}_w(z)$ and
$\hat{d}_v(z+1)=\tilde{v}(-z-1)\hat{d}_v(z)$,
$\hat{d}_w(z+1)=\tilde{w}(-z-1)\hat{d}_w(z)$
it follows that
\begin{align}\label{ado1}
\frac{\tilde{\Delta}_+(x+e_{\{ 1,\ldots ,r\} })}{\tilde{\Delta}_+(x)} &=
\tilde{V}_{ \{ 1,\ldots ,r\} ,\, \{ r+1,\ldots ,n\} } (x) \\
\intertext{and}
\label{ado2}
\frac{\hat{\Delta}_+(x+e_{\{ 1,\ldots ,r\} })}{\hat{\Delta}_+(x)} &=
\tilde{V}_{ \{ 1,\ldots ,r\} ,\, \{ r+1,\ldots ,n\} }(-x-e_{\{ 1,\ldots ,r\}
}),   \end{align}
where $\tilde{V}_{J,\, K}(x)$ is the same in Eq. \eqref{dualeqs} (with all
signs $\varepsilon_j$, $j\in J$ being positive).

\begin{thm}[normalization~1: duality]\label{norm1}
One has
\begin{equation}\label{norm1eq}
 p_\lambda (\beta \rho^*) =
\frac{\tilde{\Delta}_+(\rho+\lambda)}{\tilde{\Delta}_+ (\rho)} .
\end{equation}
\end{thm}
\begin{pf}
Expanding the l.h.s. and the r.h.s. of Recurrence relation~\eqref{recr}
in monomial symmetric functions and comparing the coefficients of the leading
monomial $m_{\lambda +e_{ \{ 1,\ldots ,r\}}}(x)$ at both sides of the equation
leads to a relation between $p_{\lambda +e_{\{ 1,\ldots ,r\} }} (\beta \rho^*)$
and $p_{\lambda }(\beta \rho^*)$:
\begin{equation}
 1/p_{\lambda }(\beta \rho^*) =
\tilde{V}_{ \{ 1,\ldots ,r\} ,\, \{ r+1,\ldots ,n\} } (\rho +\lambda ) /
p_{\lambda +e_{\{ 1,\ldots ,r\} }}(\beta \rho^*) .
\end{equation}
One can rewrite this relation with the aid of Eq.~\eqref{ado1} in the form
\begin{equation}\label{quo}
\frac{ p_\lambda (\beta \rho^*) }{ \tilde{\Delta}_+ (\rho +\lambda )}=
\frac{ p_{\lambda +e_{\{ 1,\ldots ,r\} }} (\beta \rho^*) }
     { \tilde{\Delta}_+ (\rho +\lambda +e_{\{ 1,\ldots ,r\} })} .
\end{equation}
Because any vector in the cone $\Lambda$ \eqref{monomial} can be written as a
nonnegative integral combination of the vectors $e_{\{ 1,\ldots ,r\}}$,
$r=1,\ldots ,n$,
it follows from (iterated application of) Eq.~\eqref{quo} that the quotient
$p_\lambda (\beta \rho^*)/\tilde{\Delta}_+ (\rho +\lambda )$ does not depend on
the choice of $\lambda \in \Lambda$. Evaluation of this quotient in $\lambda
=0$ (so $p_\lambda =1$) then entails Eq.~\eqref{norm1eq}.
\end{pf}
Theorem~\ref{norm1} manifestly demonstrates that $p_\lambda (\beta \rho^*)$ is
indeed nonzero for generic parameters (Lemma~\ref{nonzero}). One observes in
particular that $p_\lambda (\beta \rho^*)$ is positive
for nonnegative values of the parameters $g, g_r$
(r=0,1,2,3). (For positive parameters this is immediate from
the expression in Eq.~\eqref{norm1eq},
whereas for $g$ or $g_r$ equal to zero one needs to consider the
limit $g\rightarrow 0_+$ or $g_r\rightarrow 0_+$, respectively).

Let us next turn to the computation of the (ortho)normalization constants
$\langle p_\lambda ,p_\lambda \rangle_\Delta$ (for nonnegative values of the
parameters $g,g_0,\ldots ,g_3$). In the case $\lambda =0$ (so $p_\lambda =1$),
the corresponding (generalized Selberg-type) integral was evaluated by
Gustafson in Ref.~\cite{gus:generalization}. In our notation, Gustafson's
result may be rephrased as (cf. Remark~{\em ii.}, below)
\begin{equation}\label{gust}
\langle 1, 1 \rangle_\Delta =
 2^n n!\, \tilde{\Delta}_+ (\rho)\, \hat{\Delta}_+ (\rho) .
\end{equation}
The following theorem generalizes this (constant term) formula to the case of
Koorn\-win\-der-Macdonald polynomials of arbitrary degree.

\begin{thm}[normalization~2: orthonormality]\label{norm2}
Let the parameters $g, g_r$, $r=0,\ldots ,3$ be nonnegative. Then one has
\begin{equation}\label{ipid}
\langle p_\lambda ,p_\lambda \rangle_\Delta = 2^n n!\,
\tilde{\Delta}_+ (\rho +\lambda ) \hat{\Delta}_+ (\rho +\lambda ).
\end{equation}
\end{thm}
\begin{pf}
Set $\tilde{\Delta}(x)=\tilde{\Delta}_+(x)\, \tilde{\Delta}_+(-x)$
with $g^*,g_r^*$ replaced by $g,g_r$
(so $\tilde{\Delta} (x)$ equals $\Delta (\beta x)$, cf. Eq.~\eqref{weight})
and let
$\langle \cdot ,\cdot \rangle_{\tilde{\Delta}}$ be the inner product determined
by (cf. Eq.~\eqref{ip})
\begin{equation}\label{ip2}
\langle \tilde{m}_\lambda , \tilde{m}_{\lambda^\prime}
                    \rangle_{\tilde{\Delta}}   =
\left( \frac{\alpha\beta}{2\pi} \right)^n
\int_{-\frac{\pi}{\alpha\beta}}^{\frac{\pi}{\alpha\beta}}\cdots
\int_{-\frac{\pi}{\alpha\beta}}^{\frac{\pi}{\alpha\beta}}
\tilde{m}_\lambda (ix) \, \overline{\tilde{m}_{\lambda^\prime}(ix)}\,
\tilde{\Delta} (ix)\,  dx_1\cdots dx_n
\end{equation}
(where $\tilde{m}_\lambda (x)=m_\lambda (\beta x)/m_\lambda (0)$, cf.
Definition~\eqref{reno}).
Clearly, the polynomials $\tilde{p}_\lambda$ \eqref{reno} are orthogonal with
respect to $\langle \cdot ,\cdot \rangle_{\tilde{\Delta}}$ (cf.
Eq.~\eqref{ort}). Furthermore, one has (using
Definition~\eqref{reno} and Theorem~\ref{norm1})
\begin{equation}\label{nl}
N_\lambda \equiv \langle \tilde{p}_\lambda ,\tilde{p}_\lambda
\rangle_{\tilde{\Delta}}
= \frac{1}{| p_\lambda (\beta\rho^*)|^2}
\langle p_\lambda ,p_\lambda \rangle_\Delta =
\left(
\frac{\tilde{\Delta}_+(\rho )}{\tilde{\Delta}_+(\rho +\lambda )}\right) ^2\:
\langle p_\lambda ,p_\lambda \rangle_\Delta .
\end{equation}

By evaluating both sides of the identity
\begin{equation}
\langle E^*_r \tilde{p}_\lambda ,
        \tilde{p}_{\lambda +e_{\{ 1,\ldots ,r\}}} \rangle_{\tilde{\Delta}} =
\langle \tilde{p}_\lambda , E^*_r
        \tilde{p}_{\lambda +e_{\{ 1,\ldots ,r\}}} \rangle_{\tilde{\Delta}}
\end{equation}
using Recurrence relation~\eqref{recr} and the orthogonality of the polynomials
$\tilde{p}_\lambda$, $\lambda \in \Lambda$, we arrive at the following relation
between
$N_\lambda$ ($ =\langle \tilde{p}_\lambda ,\tilde{p}_{\lambda}
\rangle_{\tilde{\Delta}}$)
and
$N_{\lambda +e_{\{ 1,\ldots ,r\}}}$:
\begin{equation}
V_{\{1,\ldots ,r\} ,\{r+1,\ldots ,n\} }(\rho +\lambda )
N_{\lambda +e_{\{ 1,\ldots ,r\}}}  =
V_{\{1,\ldots ,r\} ,\{r+1,\ldots ,n\} }
(-\rho -\lambda -e_{\{ 1,\ldots ,r\}}) N_\lambda  .
\end{equation}
With the aid of Eqs. \eqref{ado1} and \eqref{ado2} one rewrites
this relation in the form
\begin{equation}
N_\lambda
\frac{\tilde{\Delta}_+(\rho +\lambda )}
     {\hat{\Delta}_+(\rho +\lambda )}=
N_{\lambda +e_{\{ 1,\ldots ,r\}} }
\frac{\tilde{\Delta}_+(\rho +\lambda +e_{\{ 1,\ldots ,r\}} )}
     {\hat{\Delta}_+(\rho +\lambda +e_{\{ 1,\ldots ,r\}} )},
\end{equation}
which after invoking of Eq.~\eqref{nl} goes over in
\begin{equation}\label{iprel}
\frac{\langle p_\lambda ,p_\lambda \rangle_\Delta}
     {\tilde{\Delta}_+(\rho +\lambda )\, \hat{\Delta}_+(\rho +\lambda )}=
 \frac{\langle p_{\lambda +e_{\{ 1,\ldots ,r\}}} ,
               p_{\lambda +e_{\{ 1,\ldots ,r\}}} \rangle_\Delta}
     {\tilde{\Delta}_+(\rho +\lambda +e_{\{ 1,\ldots ,r\}} )\,
     \hat{\Delta}_+(\rho +\lambda +e_{\{ 1,\ldots ,r\}} )} .
\end{equation}
It follows from Eq.~\eqref{iprel} that the quotient
$\langle p_\lambda ,p_\lambda \rangle_\Delta /
(\tilde{\Delta}_+(\rho +\lambda )\, \hat{\Delta}_+(\rho +\lambda ) )$
does not depend on $\lambda \in \Lambda$~\eqref{monomial}
(cf. the proof of Theorem~\ref{norm1}).
Combining this with Gustafson's result \eqref{gust} for $\lambda =0$ yields
Eq.~\eqref{ipid}.
\end{pf}

{\em Remarks:} {\em i.} Both the evaluations for $p_\lambda (\beta \rho^*)$
(Theorem~\ref{norm1}) and $\langle p_\lambda ,p_\lambda \rangle_\Delta$
(Theorem~\ref{norm2}) were conjectured by Macdonald in Ref.~\cite{mac:some}
(without imposing the Self-duality condition~\eqref{self-dual}).
The formula in Theorem~\ref{norm1} is sometimes referred to as the evaluation
or specialization formula/conjecture, and the formula in
Theorem~\ref{norm2} is also known as the inner product
identity/conjecture. In the special case that $\lambda =0$, the latter formula
is also called the constant term formula (as it amounts to explicit computation
of the constant term in the Fourier decomposition of $\Delta (ix)$).

{\em ii.} Equation~\eqref{gust} boils down to an explicit evaluation of the
generalized Selberg-type integral
\begin{equation}\label{selb}
\left( \frac{\alpha}{2\pi} \right)^n
\int_{-\pi/\alpha}^{\pi/\alpha}\!\!\!\!\!\!\cdots
               \int_{-\pi/\alpha}^{\pi/\alpha}
 \Delta (ix)\,  dx_1\cdots dx_n \; = \; 2^n n!\, \tilde{\Delta}_+ (\rho)\,
\hat{\Delta}_+ (\rho) .
\end{equation}
By cancelling common factors in the numerator and denominator, it is not
difficult to rewrite the r.h.s. of Eq.~\eqref{selb} as
\begin{equation}
2^n n!\, \prod_{1\leq j\leq n}
\frac{ (t;q)_\infty (a_0 a_1 a_2 a_3 t^{n+j-2};q)_\infty }
     { (q;q)_\infty (t^j;q)_\infty
       \prod_{0\leq r < r^\prime \leq 3} (a_r a_{r^\prime}t^{j-1};q)_\infty },
\end{equation}
where
$t=q^g$, $a_0=q^{g_0}$, $a_1=-q^{g_1}$, $a_2=q^{g_2+1/2}$,
and $a_3=-q^{g_3+1/2}$ (cf. Eq.~\eqref{parel}).
This is the expression for the evaluation constant
found by Gustafson~\cite{gus:generalization}.

\section{Proof of the duality theorem}\label{sec6}
In this section we will prove the Duality relations~\eqref{dr}
(Theorem~\ref{dthm}) by performing induction on $\mu$.
It is immediate from the definition of the (renormalized)
Koornwinder-Macdonald polynomials that the duality relations hold for all
$\lambda \in \Lambda$~\eqref{monomial}
if $\mu=0$, as in that case both sides of Eq.~\eqref{dr}
are identical to one. Now, let $\omega\in \Lambda$ be nonzero and let us assume
as induction hypothesis that Eq~\eqref{dr} is valid for all
$\lambda\in \Lambda$ and all $\mu\in \Lambda$ with $\mu < \omega$.
We shall prove that this implies that Eq~\eqref{dr} also holds
for $\mu = \omega$ (and all $\lambda \in \Lambda$).

Since $\omega\in \Lambda$ is nonzero, there must exist an
$s\in \{ 1,\ldots ,n-1\}$ such that $\omega_s >\omega_{s+1}$ ($\geq 0$).
Hence, $\nu \equiv \omega -e_{\{ 1,\ldots ,s\}}\in \Lambda$ \eqref{monomial}
(recall $e_{\{ 1,\ldots ,s\}} \equiv e_1+\cdots + e_s$).
The most important step in the duality proof consists of demonstrating (of
course without relying on Theorem~\ref{dthm} or Theorem~\ref{recthm}) that
$\tilde{p}_\nu (x) $ satisfies the Recurrence relation \eqref{recr} for $r=s$.
The proof of this statement mimics the derivation
presented in Section~\ref{sec4} of the recurrence relations
(Theorem~\ref{recthm}) starting from the duality relations
(Theorem~\ref{dthm}).
At some points, however, it is necessary to
adapt the arguments given there since the starting point now
is the above induction hypothesis rather than the duality theorem.

As before (cf. Section~\ref{sec4}), substituting $x=\alpha (\rho + \nu )$ in
the $s$-th difference equation for $p^*_\lambda (x)$ entails (after
renormalizing and restricting the sum with the aid of Lemma~\ref{res} to those
index sets and sign configurations such that $\nu +e_{\varepsilon J}\in
\Lambda$)
\begin{equation}\label{dfe}
E^*_s  ( \rho^* +\lambda )\: \tilde{p}^*_\lambda (\rho +\nu )=\!\!\!
\sum\begin{Sb} 0\leq |J|\leq s,\;  \varepsilon_j=\pm 1,\, j\in J \\
                \nu +e_{\varepsilon J} \in \Lambda  \end{Sb}\!\!\!\!\!\!
\tilde{U}_{J^c,\, s-|J|}(\rho +\nu)\,
\tilde{V}_{\varepsilon J,\, J^c}(\rho+\nu)\,
\tilde{p}^*_\lambda (\rho+\nu+ e_{\varepsilon J}) .
\end{equation}
Let us assume for the moment that $\lambda \leq \omega$. Then we may use the
induction hypothesis to rewrite Eq.~\eqref{dfe} as
\begin{equation}\label{re}
E^*_s  ( \rho^* +\lambda )\: \tilde{p}_\nu (\rho^* +\lambda )=\!\!\!
\sum\begin{Sb} 0\leq |J|\leq s,\;  \varepsilon_j=\pm 1,\, j\in J \\
                \nu +e_{\varepsilon J} \in \Lambda  \end{Sb}\!\!\!\!\!\!
\tilde{U}_{J^c,\, s-|J|}(\rho +\nu)\,
\tilde{V}_{\varepsilon J,\, J^c}(\rho+\nu)\,
\tilde{p}_{\nu+ e_{\varepsilon J}}(\rho^*+\lambda ) .
\end{equation}
To obtain Eq.~\eqref{re} from Eq.~\eqref{dfe},
we have used for the l.h.s. that
$ \tilde{p}^*_\lambda (\rho +\nu ) =  \tilde{p}_\nu (\rho^* +\lambda )$
because $\nu (\equiv \omega -e_{ \{ 1,\ldots ,s\}}) < \omega$, and for the
r.h.s. that
$\tilde{p}^*_\lambda (\rho+\nu+ e_{\varepsilon J})=
\tilde{p}_{\nu+ e_{\varepsilon J}}(\rho^*+\lambda )$ because either
1. $\lambda < \omega$, or 2. $\nu +e_{\varepsilon J}<\omega$ (this happens when
$e_{\varepsilon J}\neq e_{ \{ 1,\ldots ,s\}}$), or
3. $\lambda =\nu +e_{\varepsilon J}= \omega$ (this happens when $\lambda
=\omega$
and $e_{\varepsilon J}= e_{ \{ 1,\ldots ,s\}}$).
In the third case one has that
$\tilde{p}_\omega^*(\rho +\omega)= \tilde{p}_\omega(\rho^*+\omega)$ trivially,
because of the Self-duality assumption~\eqref{self-dual}
(which implies that $\tilde{p}_\omega^*(x)= \tilde{p}_\omega(x)$ and
$\rho^*=\rho$).
It is precisely at this point (and only at this point) that we actually use the
Self-duality condition (cf. Remark~\ref{sub2} of Section~\ref{sec7}).

So far we have shown that Eq.~\eqref{re} holds for all $\rho^*+\lambda$ with
$\lambda \leq \omega$. To see that equality actually holds for all values of
the argument (and thus proving the $s$-th recurrence relation for
$\tilde{p}_\nu (x)$),
we need two lemmas.
\begin{lem}\label{lem1} One has
\begin{equation}\label{expansion}
E^*_s  ( x )\: \tilde{p}_\nu (x )=
\sum_{\mu\in\Lambda ,\, \mu \leq \omega} c_\mu\: \tilde{p}_\mu (x )
\;\;\;\;\;\;\; \text{with}
\;\;\;\;\;\;\; c_\mu\in {\Bbb C} .
\end{equation}
\end{lem}
\begin{pf}
It is clear from the definitions that $E_s^*(x)\: \tilde{p}_\nu(x)$ is an
even and permutation invariant trigonometric polynomial.
Hence, it can
be expanded as a finite linear combination of monomials
$\tilde{m}_\lambda(x)$, $\lambda\in \Lambda$ (where $\tilde{m}_\lambda(x)\equiv
m_\lambda({\beta x})/m_\lambda (0)$, cf. Definition~\eqref{reno}).
It follows from the asymptotics at infinity
that in this expansion only monomials
$\tilde{m}_\lambda(x)$ with $\lambda \leq \omega$ occur with a nonzero
coefficient.
To see this we set $x=Ry$
with $y_1>y_2>\cdots >y_n>0$ and notice that
\begin{gather}\label{monas}
\tilde{m}_\lambda (Ry) \sim  e^{\alpha\beta R \langle \lambda , y \rangle }
\;\;\;\;\;\;\;\; \text{for}\;\;\; R\rightarrow \infty \\
\intertext{(where
$\langle \lambda ,y\rangle \equiv\sum_{1\leq j\leq n} \lambda_j y_j$
and $\sim$ denotes proportionality). One furthermore has}
E^*_s(Ry)\sim e^{\alpha\beta R \langle e_{\{ 1,\ldots ,s\} },y \rangle },
\;\;\;\;\;\;\;\;\;
\tilde{p}_\nu(Ry) \sim  e^{\alpha\beta R \langle \nu , y \rangle}
\;\;\;\;\;\;\;\; \text{for}\;\;\; R\rightarrow \infty .
\end{gather}
In obtaining the asymptotics for $\tilde{p}_\nu(Ry)$ we have used
the fact that
\begin{equation}\label{smaller}
\lambda^\prime  \leq \lambda \;\;\;\;\; \text{iff}\;\;\;\;\;
 \langle \lambda^\prime , y \rangle \leq \langle \lambda , y \rangle
\;\;\;
\forall y\in {\Bbb R}^n\;\;\text{with}\;\; y_1>\cdots>y_n>0 .
\end{equation}
By comparing the asymptotics of
$E_s^*(Ry)\: \tilde{p}_\nu(Ry)\sim
\exp (\alpha\beta R \langle \omega , y \rangle )$ (recall
$\omega =\nu +e_{\{ 1,\ldots ,s\}}$) with that of
$\tilde{m}_\lambda (Ry)$~\eqref{monas},
one infers---using Eq.~\eqref{smaller}---that
\begin{equation}\label{monex}
E^*_s  ( x )\: \tilde{p}_\nu (x )=
\sum_{\mu\in\Lambda ,\, \mu \leq \omega} \hat{c}_\mu\: \tilde{m}_\mu (x )
\;\;\;\;\;\;\; \text{with}
\;\;\;\;\;\;\; \hat{c}_\mu \in {\Bbb C} .
\end{equation}
The lemma is now immediate from Eq.~\eqref{monex} and the fact that the
transition matrix relating the bases
$\{\tilde{p}_\lambda(x)\}_{\lambda \in \Lambda}$ and
$\{\tilde{m}_\lambda(x)\}_{\lambda \in \Lambda}$ is upper triangular
($\tilde{p}_\lambda(x) =
\sum_{\lambda^\prime\in\Lambda ,\, \lambda^\prime \leq \lambda}
\tilde{c}_{\lambda ,\lambda^\prime}
\,\tilde{m}_{\lambda^\prime}(x)$ and
$\tilde{m}_\lambda(x) =
\sum_{\lambda^\prime \in \Lambda ,\,\lambda^\prime \leq \lambda}
\tilde{c}^{inv}_{\lambda ,\lambda^\prime}
\,\tilde{p}_{\lambda^\prime}(x)$).
\end{pf}

\begin{lem}\label{lem2} For generic parameters $g_,g_0,\ldots ,g_3$
and scale factors $\alpha ,\beta$,
the matrix $[\tilde{p}_\mu (\rho^*+\lambda )]_{\mu,\lambda}$,
with $\mu ,\lambda \in \Lambda$ and $\mu ,\lambda \leq \omega$,
is invertible.
\end{lem}
\begin{pf}
{}From the asymptotics
$\tilde{m}_\mu (\rho^*+\lambda )\sim
\exp (\alpha\beta \langle \mu , \rho^*+\lambda \rangle)=
q^{-\langle \mu ,\rho^*+\lambda \rangle}$ for $g^*(=g)\rightarrow\infty$
and the fact that
$\text{det}[q^{-\langle \mu ,\lambda \rangle}]_{\lambda ,\mu\leq\omega}\neq 0$
for generic $q$ (see Remark at the end of this section),
it follows that the determinant
$\text{det}[\tilde{m}_\mu (\rho^*+\lambda )]_{\mu,\lambda\leq\omega}$
is not identically zero. The same is then true for
$\text{det}[\tilde{p}_\mu (\rho^*+\lambda )]_{\mu,\lambda\leq\omega}$,
as the transition matrix relating the bases
$\{\tilde{p}_\lambda(x)\}_{\lambda \in \Lambda}$ and
$\{\tilde{m}_\lambda(x)\}_{\lambda \in \Lambda}$
(and therefore also the matrix relating
$\text{det}[\tilde{p}_\mu (\rho^*+\lambda )]_{\mu,\lambda\leq\omega}$
and
$\text{det}[\tilde{m}_\mu (\rho^*+\lambda )]_{\mu,\lambda\leq\omega}$)
is (upper) triangular
(cf. the proof of Lemma~\ref{lem1}).
Thus, in view of the analyticity in the parameters and in the scale factors
$\alpha ,\beta$, it is clear that
$\text{det}[\tilde{p}_\mu (\rho^*+\lambda )]_{\mu,\lambda\leq\omega}\neq 0$
generically, which proves the Lemma.
\end{pf}
With these two lemmas we are finally in the position
to show that Eq.~\eqref{re} not only holds for $\rho^*+\lambda$,
$\lambda \leq \omega$ but in fact for all values of the argument.
Substituting $x=\rho^*+\lambda$ in Expansion~\eqref{expansion} and
comparing with Eq.~\eqref{re}, tells us---using Lemma~\ref{lem2}---that
the expansion coefficient $c_\mu$ equals
$\tilde{U}_{J^c,\, s-|J|}(\rho +\nu)\,
\tilde{V}_{\varepsilon J,\, J^c}(\rho+\nu) $ if $\mu =\nu +e_{\varepsilon J}$
(for some index set $J$ with $0\leq |J|\leq s$ and configuration of signs
$\varepsilon_j$, $j\in J$), and that $c_\mu =0$ otherwise.
It is clear that with such expansion coefficients $c_\mu$,
Eq.~\eqref{expansion} goes over in
the $s$-th recurrence relation for $\tilde{p}_\nu (x )$:
\begin{equation}\label{recs}
E^*_s ( x )\: \tilde{p}_\nu (x )=\!\!\!
\sum\begin{Sb} 0\leq |J|\leq s,\;  \varepsilon_j=\pm 1,\, j\in J \\
                \nu +e_{\varepsilon J} \in \Lambda  \end{Sb}\!\!\!\!\!\!\!
\tilde{U}_{J^c,\, s-|J|}(\rho +\nu)\,
\tilde{V}_{\varepsilon J,\, J^c}(\rho+\nu)\,
\tilde{p}_{\nu+ e_{\varepsilon J}}(x ) .
\end{equation}

It follows in particular from the Recurrence relation~\eqref{recs} that
Eq.~\eqref{re} is valid for {\em all} $\lambda \in \Lambda$
(and not just for $\lambda \leq \omega$).
Furthermore, if we subtract Eq.~\eqref{re} from Eq.~\eqref{dfe}, then all
terms in the difference---except the leading term in the r.h.s.
corresponding to the index set $J=\{ 1,\ldots ,s\}$ with all signs
$\varepsilon_j$
($j\in J$)
positive---manifestly cancel each other in view of the
induction hypothesis
(recall $\tilde{p}^*_\lambda (\rho +\nu )=\tilde{p}_\nu (\rho^* +\lambda )$
because $\nu < \nu +e_{\{1,\ldots ,s\}}= \omega$, and
$\tilde{p}^*_\lambda (\rho+\nu+ e_{\varepsilon J})=
\tilde{p}_{\nu+ e_{\varepsilon J}}(\rho^*+\lambda )$ for
$e_{\varepsilon J}\neq e_{\{ 1,\ldots ,s\}}$ because
$\nu +e_{\varepsilon J}<\nu +e_{\{ 1,\ldots ,s\}}=\omega$
if $e_{\varepsilon J}\neq e_{\{ 1,\ldots ,s\}}$).
Hence, we end up with the equation (recall $U_{K,\, 0}=1$)
\begin{equation}
0=V_{ \{ 1,\ldots ,s\},\; \{s+1,\ldots ,n\} } (\rho +\nu )\;
\left( \tilde{p}_\lambda^*(\rho+\omega)-\tilde{p}_\omega (\rho^*+\lambda)
\right) .
\end{equation}
For generic parameters
$V_{ \{ 1,\ldots ,s\},\; \{s+1,\ldots ,n\} } (\rho +\nu )$ is nonzero
(Lemma~\ref{res}), so we find that
$\tilde{p}_\lambda^*(\rho+\omega)=\tilde{p}_\omega (\rho^*+\lambda)$
for all $\lambda \in\Lambda$
(first for generic parameters and then for all parameters in view of
the analyticity in the parameters).

This completes the induction step and thus the proof of Theorem~\ref{dthm}.

{\em Remark:} In the proof of Lemma~\ref{lem2} we have used that
$\text{det}[q^{-\langle \mu ,\lambda \rangle}]_{\lambda ,\mu\leq\omega}\neq 0$
for generic $q$. This follows from a more general result due to Koornwinder
\cite{koo:self} (see also Ref.~\cite[Ch. 6]{mac:symmetric})
stating that for distinct vectors $v_1,\ldots ,v_N\in {\Bbb R}^n$
(and generic $q$)
\begin{equation}\label{det}
\text{det} [q^{\langle v_j, v_k\rangle} ]_{1\leq j,k\leq N} =
\sum_{\sigma \in S_N} (-1)^\sigma
q^{\langle v_j, v_{\sigma (j)}\rangle} \neq 0.
\end{equation}
To see that the terms $q^{\langle v_j, v_{\sigma (j)}\rangle}$ in
Eq.~\eqref{det} cannot all cancel each other one uses the estimate
\begin{equation}
\sum_{1\leq j\leq N} \langle v_j, v_{\sigma (j)}\rangle \leq
\sum_{1\leq j\leq N} \langle v_j, v_{j}\rangle
\end{equation}
(immediate from the Cauchy-Schwarz inequality (twice)),
with equality holding only when $\sigma =\text{id}$.
It follows from this estimate,
first for $q\rightarrow \infty$ and hence for generic (say positive) $q$
by analyticity, that the determinant in Eq.~\eqref{det} is indeed nonzero.

\section{Concluding remarks}\label{sec7}
\subsection{Structure of the difference equations}
It has already been pointed out in Remark~{\em i.} of Section~\ref{sec3} that
the structure of the Difference equations \eqref{diffeqr} is that of a system
of eigenvalue equations for a family of commuting difference operators
$D_r$~\eqref{ados} ($r=1,\ldots ,n)$. The heighest-order part of $D_r$ is of
the form (recall $U_{K,\, 0}=1$)
\begin{equation}
\sum\begin{Sb} J\subset \{ 1,\ldots ,n\} ,\, |J|= r \\
               \varepsilon_j=\pm 1,\; j\in J      \end{Sb}\!\!\!\!
V_{\varepsilon J,\, J^c}(x)\,
T_{\varepsilon J,\, \beta}.
\end{equation}
The coefficients $V_{\varepsilon J,\, J^c}(x)$ \eqref{diffeqr} are related to
the weight function $\Delta(x)$ \eqref{weight} via (cf. Eq.~\eqref{ado1})
\begin{equation}
V_{ \{ 1,\ldots ,r\} ,\, \{ r+1,\ldots ,n\} } (x)=
\Delta_+(x+\beta e_{\{ 1,\ldots ,r\} })/\Delta_+(x)  ,
\end{equation}
where $\Delta_+(x)=\tilde{\Delta}_+(\beta^{-1}x)$ with $g^*, g_r^*$
replaced by $g,g_r$ (so $\Delta (x)=\Delta_+(x)\Delta_+(-x)$).
The lower-order terms of $D_r$ \eqref{ados} consist of a similar sum over the
index sets $J\subset \{1,\ldots ,n\}$
with cardinality smaller than $r$ and with each term
$V_{\varepsilon J,\, J^c}(x)\, T_{\varepsilon J,\, \beta}$
being multiplied by the function $U_{J^c,\, r-|J|}(x)$.

Now, for $\lambda =0$ the r.h.s. of the Difference equation \eqref{diffeqr}
vanishes because $E_r(\rho)=0$ \cite{die:commuting}. Equivalently, one could
say that the difference operators $D_1,\ldots ,D_n$ \eqref{ados} annihilate
constant functions. (For $D=D_1$ this is clear from Eq.~\eqref{dop1}.)
This property of the difference operators/equations
gives rise to the following set of relations between functions
$V_{\varepsilon J,\, I}$ and $U_{K,\, p}$:
\begin{equation}\label{UVrel}
\sum\begin{Sb} J\subset \{ 1,\ldots ,n\} ,\, 0\leq |J|\leq r \\
               \varepsilon_j=\pm 1,\; j\in J      \end{Sb}\!\!\!\!
U_{J^c,\, r-|J|}\,  V_{\varepsilon J,\, J^c} =0, \;\;\;\;\;\;
r=1,\ldots ,n.
\end{equation}
These functional relations (together with the condition that
$U_{K,\, p}=1$ for $p=0$) actually determine $U_{K,\, p}$
completely in terms of $V_{\varepsilon J,\, I}$.
For instance, for $r=1$ Eq.~\eqref{UVrel} entails
\begin{equation*}
\begin{align*}
U_{K,\, 1} =& - \sum_{k\in K}
( V_{\{ k\} ,\, K\setminus \{ k\} } +
 V_{- \{ k\} ,\, K\setminus \{ k\} } )  \\
 = &
-\sum_{k\in K} \Bigl(
 w(x_k) \prod_{k^\prime \in K,\, k^\prime \neq k}
v(x_k+x_{k^\prime})\, v(x_k-x_{k^\prime})\; + \\
& \;\;\;\;\;\;\;\; w(-x_k) \prod_{k^\prime \in K,\, k^\prime \neq k}
v(-x_k+x_{k^\prime})\, v(-x_k-x_{k^\prime}) \Bigr) .
\end{align*}
\end{equation*}
(Eq.~\eqref{UVrel} immediately yields $U_{K,\, 1}$ for $K=\{ 1,\ldots ,n\}$;
the case of general $K$ then follows by renumbering.)

More generally, one obtains $U_{K,p}$ for general $p$ from Eq.~\eqref{UVrel}
with $r=p$ by performing induction on $p$. This way one finds
\begin{equation*} U_{K,p}\; =
\sum\begin{Sb}
    \emptyset\subsetneq L_1\subsetneq \cdots \subsetneq L_m\subset K,\;
     1\leq m\leq p      \\
     |L_m|=p,\; \varepsilon_l=\pm 1,\, l\in L_m
    \end{Sb} (-1)^m\;\;
V_{\varepsilon L_1 ,\, K\setminus L_1}\,
V_{\varepsilon (L_2 \setminus L_1),\, K\setminus L_2}\cdots
V_{\varepsilon (L_m \setminus L_{m-1}),\, K\setminus L_m} .
\end{equation*}
It turns out that this expression for $U_{K,\, p}$ can be rewritten in
the more compact form that was used in Eq.~\eqref{diffeqr}
\cite{die:difference,die:diagonalization}.
The equality of both expressions for
$U_{K,p}$ hinges on a system of functional equations for $v(z)$
satisfied by $v(z) =\text{sh}\, \frac{\alpha}{2} (\beta g +z)/
\text{sh} \, (\frac{\alpha}{2} z)$.

\subsection{Dropping the self-duality condition}\label{sub2}
In Section~\ref{sec6} we provided a proof of the duality relations
for the Koornwinder-Macdonald polynomials
(Theorem~\ref{dthm}) with parameters subject to the Self-duality
condition~\eqref{self-dual}.
This condition effectively reduces the number of independent parameters from
five to four (not counting the scale factors $\alpha$ and $\beta$).
It is expected (and conjectured by Macdonald~\cite{mac:some}), however, that
the Duality relations~\eqref{dr} are true for the full five-parameter
Koornwinder-Macdonald family (thus generalizing the state of affairs
for $n=1$).
Should one succeed in proving the Duality theorem~\ref{dthm} without
restrictions on the parameters, then automatically all other results of
Sections \ref{sec4} and \ref{sec5} carry over to this
slightly more general situation
(with the proofs given applying verbatim).

A careful examination of the proof of Theorem~\ref{dthm}
given in Section~\ref{sec6}
reveals that the only step requiring invoking of
the Self-duality condition~\eqref{self-dual} has been the derivation of
Eq.~\eqref{re} from Eq.~\eqref{dfe} (with the aid of the induction hypothesis).
At that point we needed that
$\tilde{p}_\omega^*(\rho +\omega) =\tilde{p}_\omega (\rho^*+\omega)$, which is
trivial for self-dual polynomials
(because then
$\tilde{p}_\omega^*(x) =\tilde{p}_\omega (x)$ and $\rho^*=\rho$), but which
requires a proof when the self-duality condition is dropped.
If one would be able to prove the relation
$\tilde{p}_\omega^*(\rho +\omega) =\tilde{p}_\omega (\rho^*+\omega)$
($\omega\in \Lambda$) for arbitrary parameters, then Theorem~\ref{dthm}
(and thus all other theorems in Sections~\ref{sec4} and \ref{sec5})
would follow immediately for the complete five-parameter Koornwinder-Macdonald
family.

There might be an alternative approach. If one chooses to {\em define} the
renormalized Koornwinder-Macdonald polynomials $\tilde{p}_\lambda (x)$
as $p_\lambda (\beta x)$ divided by the constant
in the r.h.s. of Eq.~\eqref{norm1eq},
then the derivation of Eq.~\eqref{recs}
can be established starting from Eq.~\eqref{dfe} by assuming
$\lambda < \omega$ and applying the induction hypothesis to arrive at
Eq.~\eqref{re} (because $\lambda <\omega$ it is now not necessary to
invoke the self-duality condition).
Next we bring all terms in the r.h.s. of Eq.~\eqref{re}
with $|J|=s$ to the l.h.s. and
employ a version of Lemma~\ref{lem1} and Lemma~\ref{lem2}
with strict inequalities ($\mu ,\lambda< \omega$)
to arrive at Eq.~\eqref{recs}. (By bringing the terms with $|J|=s$
to the other side one ensures that the resulting function
in the l.h.s. can be expanded in
polynomials $\tilde{p}_\mu (x)$ with $\mu <\omega$.)
This proves the induction step.
However, to check now that the duality relations hold for $\mu =0$
amounts to proving Theorem~\ref{norm1}.
Hence, the upshot is that all results of the paper can be extended to
arbitrary parameters once Eq.~\eqref{norm1eq}
(Theorem~\ref{norm1}) has been verified.
Or, in other words,
the Duality relations~\eqref{dr} and the Evaluation (or specialization)
formula~\eqref{norm1eq} follow from each other.

\subsection{Affine Hecke algebras}
In Ref.~\cite{che:double}, Cherednik proved Macdonald's orthonormalization
conjectures for the Macdonald polynomials related to reduced root systems
by means of certain Hecke-algebraic techniques (leading to shift operators).
More precisely, Cherednik considered the case of admissible pairs
of the form $(R,R^\vee)$  with $R$ reduced.
Recently, these algebraic methods also resulted in a proof
of the corresponding duality relations and
renormalization formulas
(i.e., the evaluation/specialization formulas, cf. Remark~{\em i.}
of Section~\ref{sec5}.) \cite{che:macdonalds}.
For the type $A$ root system this approach reproduces
(albeit in a completely different manner) the results of
Koornwinder and Macdonald \cite{koo:self,mac:symmetric}.
For the remaining
classical reduced root systems (i.e., the types $B$, $C$, and $D$), the
$(R,R^\vee)$-type Macdonald polynomials amount to (self-dual)
Koornwinder-Macdonald polynomials with special parameters $g_r$
(type $B$: $g_0=g_2$, $g_1,g_3=0$; type $C$: $g_0=g_1$, $g_2,g_3=0$;
type $D$: $g_0,g_1,g_2,g_3=0$). (For the root systems $B_n$ and
$D_n$ one obtains only half the polynomials via the above specialization
of the parameters;
the other half of the polynomials is obtained by specializing to
parameters that are not self-dual, cf. Ref.~\cite[Section 5]{die:commuting}.)
It hence follows that for the above parameters
Refs.~\cite{che:double,che:macdonalds} provide an alternative approach
towards the proof of the Theorems~\ref{dthm}, \ref{norm1}, and \ref{norm2}.

Meanwhile, both Noumi \cite{nou:macdonald} and Macdonald \cite{mac:affine}
independently announced that Cherednik's  Hecke-algebraic techniques
may be extended to the full five-parameter Koornwinder-Macdonald family.
Hopefully this will eventually  lead to an alternative proof of the
Theorems~\ref{dthm}, \ref{norm1}, and \ref{norm2} valid for all parameters.

One aspect of the present paper seems hard to achieve via Hecke algebras,
though. It does not seem very likely that such algebraic methods will
independently reproduce our explicit formulas
for the difference equations and the recurrence relations for the
Koornwinder-Macdonald polynomials. At present Cherednik's techniques
permit concluding the existence of such difference equations and
recurrence relations (of course for special parameters)
without providing much clue as regards to their
explicit form except in the simplest cases, e.g. when $r=1$
(cf. in this connection also the statements in second paragraph of
the proof of Theorem~5.1 of Cherednik's paper Ref.~\cite{che:double}).

\subsection{Quantum groups}
Another algebraic framework in which the Koorn\-winder-Mac\-donald polynomials
appear is the representation theory of (compact) quantum groups.
Specifically, it has been demonstrated recently by Noumi and Sugitani
that for certain special values of the parameters the Koornwinder-Macdonald
polynomials may be interpreted as zonal spherical functions on compact quantum
symmetric spaces of classical type \cite{nou-sug:quantum}
(see also Ref.~\cite{nou:macdonalds} for a detailed
treatment of the type $A$ case).

In addition, it turned out that for the root system $A_{n-1}$
Macdonald's polynomials may also be
realized as vector valued characters of $U_q(sl_n)$ \cite{eti-kir:macdonalds}.
This observation has led to yet another (representation-theoretic) proof
of Koornwinder's duality and recurrence relations
for the $A$-type Macdonald polynomials \cite{eti-kir:representation}.

\subsection{Integrable systems}
It is possible to view the commuting difference operators
$D_1,\ldots ,D_n$ \eqref{ados}
as a complete set of quantum integrals for an integrable quantum mechanical
$n$-particle model \cite{die:difference,die:diagonalization}.
Similar integrable systems associated with and
diagonalized by the Macdonald polynomials related to classical root systems
are obtained via limit transitions (type $A$) or specialization
of the parameters (type $B$, $C$, $D$, and $BC$) \cite{die:commuting}.
For the type $A$ root systems the commuting difference operators
(quantum integrals) of the model were already found independently by
Ruijsenaars \cite{rui:complete} and Macdonald \cite{mac:orthogonal}.
In this special case, Ruijsenaars also studied in great detail the properties
of the corresponding classical mechanical systems \cite{rui:action}.
It is interesting to note that also at the classical level duality
relations, which were actually known even before their quantum counterparts
were discovered, play a crucial role in solving the system.

\section*{Acknowledgments} The author would like to thank M. Noumi for
explanation of his results in Ref.~\cite{nou:macdonald}
and T. H. Koornwinder for providing a copy of Macdonald's
unpublished notes \cite{mac:some}.

\bibliographystyle{amsplain}

\end{document}